\documentclass[11pt]{article}
\usepackage[utf8]{inputenc}
\usepackage{amsfonts}
\usepackage[width=17cm, left=1cm]{geometry}
\usepackage{amssymb}
\usepackage{amsmath}
\usepackage{amsthm}

\newtheorem{definition}{Definition}

\newtheorem{theorem}{Theorem}

\title{One-dimensional central measures on numberings of ordered sets}

\begin{document}

\author{A.~M.~Vershik\thanks{St.~Petersburg Department of Steklov Institute of Mathematics, St.~Petersburg State University, and Institute for Information Transmission Problems. Supported by the RSF grant 21-11-00152.}}

\maketitle

%\date{\today}

\begin{abstract}
We describe one-dimensional central measures on numberings (tableaux) of ideals of partially ordered sets (posets). As the main example, we study the poset ${\Bbb Z}_+^d$ and the graph of its finite ideals, multidimensional Young tableaux; for $d=2$, it is the ordinary Young graph. The central measures are stratified by dimension; in the paper we give a complete description of the one-dimensional stratum and prove that every ergodic central measure is uniquely determined by its frequencies. The suggested method, in particular, gives the first purely combinatorial proof of E.~Thoma's theorem for one-dimensional central measures different from the Plancherel measure (which is of dimension~$2$).
\end{abstract}

\section{Combinatorial introduction}

\subsection{General definitions. Posets, ideals, numberings}

We consider partially ordered sets (posets) $(P,\succ)$. In what follows, we assume that they are finite (or, in all nontrivial examples, countable) and, besides, finitely generated: there exists a finite set of pairwise incomparable elements that generates~$P$. For convenience, we also assume that $P$ has a~minimal element.

A (lower) ideal in a poset $P$ is a subset in $P$ such that whenever it contains an element, it contains all smaller elements.
 A {\it minimal infinite ideal} is an ideal that stops being an ideal if an arbitrary element is deleted from it. The set of all finite ideals of a poset~$P$ is a {\it distributive lattice}~$\Gamma(P)$; by Birkhoff's theorem, the converse is also true: every finitely generated distributive lattice is the lattice of ideals of a poset. The $\Bbb N$-graded graph corresponding to this lattice (the Hasse diagram of the lattice) will also be denoted by~$\Gamma(P)$.

 A  {\it numbering} of a poset $P$ (or of an ideal $I$ of $P$) is a map $\Psi:{\Bbb N}\rightarrow P$ (respectively, $\Psi:{\Bbb N}\rightarrow I$) satisfying the monotonicity condition: $x\succ y$ implies $\Psi^{-1}(x)>\Psi^{-1}(y)$; note that if
$x=\Psi(n+1)\succ \Psi(n)=y$, then $x$ is a {\it direct successor} of $y$ in the sense of the order, in other words, the interval
$(y,x)$ is empty.

The image $\Psi(\textbf{N})$ of a numbering $\Psi$, as well as $\Psi(\textbf{n})$ for $\textbf{n}=\{1,\ldots,n\}$, are
 ideals of $P$.

A numbering defines a path in the graph $\Gamma(P)$ in a natural way. By analogy with the terminology related to the Young graph, which corresponds to the poset $P={\Bbb Z}_+^2$, {\it it is natural to say that ideals of~$\Gamma(P)$ are diagrams, and numberings of ideals are tableaux corresponding to these diagrams}.

We will be interested in {\it random central numberings}, i.e., probability measures on numberings or, equivalently, measures on the path space of the graph $\Gamma(P)$; they are defined as random numberings for which the conditional probability on the finite set of paths leading to a given point $\Psi(n)$ is a~uniform measure. Thus, such a measure is a  ``maximal entropy'' measure.

\subsection{Spaces of one-dimensional ideals and one-dimensional numberings}

The set of nonempty minimal ideals of a countable poset is, in turn, a poset under the inclusion relation; it is denoted by
$\operatorname{Min}(P)$.

A {\it one-dimensional minimal ideal} in an arbitrary lattice is an ideal that, regarded as a set, is a~union of finitely many chains.

The set of one-dimensional minimal ideals with the order induced from the set $\operatorname{Min}(P)$ of all minimal ideals will be denoted by $\min_1(P)$.

Let us describe this set for our main examples.

\smallskip
a) $\operatorname{Min}({\Bbb Z}_+^2)$ is the set of all pairs of positive integers except the zero pair:
$\operatorname{Min}({\Bbb Z}_+^2)={\Bbb Z}_+^2\setminus\{(0,0)\}\bigcup \{\infty\}$; namely, a pair $(t_1,t_2)\ne (0,0)$ corresponds to the ideal consisting of the
first $t_1$ infinite rows  and the first $t_2$ infinite columns of the lattice ${\Bbb Z}_+^2$, and the symbol
${\infty}$ corresponds to the improper minimal ideal.
 In this example all proper minimal ideals are one-dimensional.

\smallskip
b) The space of one-dimensional minimal ideals of the lattice ${\Bbb Z}_+^d$ for $d>2$  has a more complicated structure.

An element of this space is indexed by a collection of $d$ Young diagrams $\nu_i$, $i=1,2\dots,  d$, of dimension $d-1$. If all diagrams except the $i$th one are empty, then this ideal is a subset of ${\Bbb Z}_+^d$ consisting of the elements whose $i$th coordinate is arbitrary (i.e., is an arbitrary nonnegative integer) and the remaining $d-1$ coordinates define a vector lying in the $i$th
$(d-1)$-dimensional Young diagram~$\nu_i$. Such one-dimensional ideals will be called
{\it irreducible}. A general one-dimensional minimal ideal in~${\Bbb Z}_+^d$ is a union of several (at most~$d$) irreducible ideals corresponding to different subsets of cardinality~$d-1$ in the set $\{1,2, \dots d\}$.

The minimal ideals of dimensions $2,3,\dots, d-1$ can be constructed in a similar way. In this paper, we do not discuss them except in the ``Comments'' section.

\smallskip
c) {\it In the case of an arbitrary poset} $P$, a minimal ideal, being a set of chains, defines a finite ideal in the finite poset of chains.
 The author does not know whether these notions have been considered in the general theory of posets or in lattice theory.

\smallskip
A numbering (tableau) of an ideal is said to be {\it one-dimensional} if either the ideal itself is one-dimensional (i.e., it is a finite ordered union of chains), or it is an infinite ordered union of chains and the numbering agrees with this ordering in the natural sense.

For example, a numbering of the entire lattice ${\Bbb Z}_+^2$ regarded as a linearly ordered union of horizontal (or vertical) chains is one-dimensional, although on the lattice itself, regarded as an improper ideal, there are non-one-dimensional numberings.

In a similar way we define the dimension of an arbitrary numbering of an ordered set (see the ``Comments'' section).

 {\it We will study random one-dimensional numberings, i.e., probability measures on the set of one-dimensional numberings. The main notion here is that of centrality of measures (see above), and the main problem is that of describing central measures on one-dimensional numberings of posets.}

\subsection{Frequencies}

Now we consider a one-dimensional numbering and define the notion of the frequency of a chain.

\begin{definition}
{\rm Let $\Psi:\textbf{N}\rightarrow P$ be a numbering of a poset $P$ and $l$ be an arbitrary chain, i.e., a linearly ordered countable subset in $P$. The {\it frequency} of $l$ is the limit
(if it exists)
$$\lim_{n\to\infty}\frac{\#\{k:k<n, \Psi(k)\in l\}}{n}=\lambda_l.$$
To a random numbering $\Psi$ there corresponds a measure $\mu_{\Lambda}$ on the space of one-dimensional numberings. For almost every individual numbering and for every chain occurring in this numbering (there are at most countably many of them), the frequency of this chain exists; the (at most countable) collection of frequencies of all chains is called the
{\it spectrum of the measure}.}
\end{definition}

It follows from general considerations of ergodic theory that for ergodic central measures the frequencies exist, i.e., every central measure determines a spectrum. Since a one-dimensional ideal is a union of chains, the finite sum of frequencies for a one-dimensional tableau must be equal to~$1$. In the case of countably many frequencies, we impose the same condition on the sum.

Now we can formulate the main question: Is every random central one-dimensional numbering uniquely determined by its spectrum? We say that a spectrum is nondegenerate if all the frequencies are distinct.

Our main result is contained in the following theorem.

 \begin{theorem}
 Every nondegenerate finite spectrum (i.e., a finite collection of distinct frequencies) determines a unique one-dimensional ergodic central measure with these frequencies.
\end{theorem}

The  nondegeneracy and finiteness conditions will be relaxed in a subsequent paper.

Before proving the theorem, we make yet another reduction of the problem for the lattices ${\Bbb Z}_+^d$, $d>1$, we are interested in.\footnote{The reduction possibly makes sense for arbitrary posets, but apparently it requires a long digression into the general theory of posets provided that such a theory exists.}

The reduction is as follows. One-dimensional minimal ideals in ${\Bbb Z}_+^d$ are unions of irreducible ideals. It is clear that the intersection of different irreducible ideals is a finite ideal. Since the frequencies of chains are preserved under finite changes, the union of the frequencies of all irreducible components coincides with the collection of frequencies of the original minimal ideal. Therefore, without loss of generality we may consider each of the irreducible minimal ideals separately.

The same consideration applies also to general one-dimensional numberings, i.e., to the case of countably many chains and, accordingly, countably many frequencies. Thus, it suffices to prove the uniqueness theorem for every central one-dimensional numbering of an irreducible ideal. After this remark, we can finally state the problem as a probability problem for finite posets. Theorem~1 immediately follows from Theorem~2 (see below), more exactly, from the positive answer to the question posed in the next section.

\section{Central measures on lattice sequences}

The following setting is not related in any way to the previous considerations and can be understood by itself.

Consider a finite poset  $P$ and the space  $P^{\infty}$ of all infinite sequences of symbols from~$P$. We say that
$\{x_n\}\in P^{\infty}$ is a {\it lattice sequence} if it satisfies the following condition: for any two elements
 $a, b \in P$ such that
$a\succ b$, and any $n\in\textbf{N}$, the inequality $$\#\{k:k<n,\, x_k=a\}>\#\{k:k<n,\, x_k=b\}$$ holds.
The set of all lattice sequences will be denoted by $R^P$.

Now we define the {\it de Finetti partitions} on the space $P^{\infty}$ and on the subset $R^P\subset P^{\infty}$.
Two sequences $\{x_n\}, \{x'_n\}$ belong to the same class of the de Finetti partition if they belong to the same orbit of the countable symmetric group
$S_{\textbf{N}}$, which acts on $P^{\infty}$ by permutations of coordinates; in other words, for some $n$
the segments  $\{x_k\}_{k=1}^n, \{x'_k\}_{k=1}^n$ have the same multisets of elements and
 $x_m=x'_m$ for $m>n$. Note that $R^P$ regarded as a subset of $P^{\infty}$ is not invariant under the action of~$S_{\textbf{N}}$, but nevertheless the de Finetti partition is well defined on~$R^P$. Denote these two de Finetti partitions by
$\xi_{P^{\infty}}$ and~$\xi_{R^P}$; the latter is the restriction of the former to~$R^P$. As we will see, this simple remark plays a crucial role in the sequel.

A central measure on $P^{\infty}$ or on $R^P$ is a (cylinder) probability measure such that for every $n$ the conditional measure on the initial segment $\{x_1,\dots, x_n\}$  given that the ``tail'' $\{x_{n+1}, x_{n+2},\dots\}$ is fixed
 (note that the set of ``beginnings'' is finite) is uniform. Thus, the notion of centrality determines a~collection of measures on an arbitrary space of sequences. Of course, the central measures themselves are different for different spaces.

An  {\it ergodic measure} on an arbitrary space of sequences is any measure  $\mu$ for which the intersection
 $\bigcap_n {\frak A}_n^{\infty}$, where ${\frak A}_n^{\infty}$ is the $\sigma$-algebra
 consisting of the sets  determined by conditions on coordinates with indices greater than $n$, is trivial.

The classical de Finetti theorem says that every ergodic central measure on $P^{\infty}$ is a product (Bernoulli) measure. Our task is to describe the ergodic central measures on $R^P$.

First of all, for every ergodic central measure on $P^{\infty}$ and, in particular, on $R^P$, we can define frequencies  $\lambda_p$ for each of the symbols $p\in P$:
      $$\lim_n \frac{\#(k<n:x_k=p)}{n}=\lambda_p\in[0,1];$$ here $\sum_{p\in P}\lambda_p=1$.

For the space of lattice sequences $R^P$, it is natural to assume that {\it the frequences agree with the order}:
$ p\succ p' \Rightarrow \lambda_p >\lambda_{p'}$.

Recall that the set of all ergodic central measures (for a graph, group, etc.) is called the {\it absolute}.

Our main technical result is as follows.

\begin{theorem}
For every system of pairwise distinct nonnegative frequencies $\Lambda=\{\lambda_p,\, p \in P\}$ that agrees with the order there exists a unique ergodic central measure $\mu_\Lambda$
on $R^P$. The set of measures $\mu_\Lambda$ exhausts the absolute.
\end{theorem}

Thus, the theorem says that the set of spectra agreeing with the order is the absolute of the space~$R^P$ of lattice sequences, and thus  there is a bijection between the Bernoulli measures and the central measures on~$R^P$.

\begin{proof}
1. First of all, note that if $P$ is finite and all the frequencies (i.e., probabilities in the Bernoulli scheme)
$\Lambda=\{\lambda_p,\, p\in P\}$ are distinct, then the set $R^P$ of lattice sequences has a positive Bernoulli measure
$\mu_{\Lambda}(R^P)$. Indeed, since $P$ is finite, it suffices to check this for a two-point poset; then we should consider the set of lattice sequences with respect to a pair of elements of~$P$ and take the intersection of these sets
 over all pairs of elements of~$P$.

Let $P$ be a two-point poset, say $P=\{1,0\}$ with $\lambda_1 > \lambda_0$. Applying Chebyshev's inequality to the random sequence $\{x_k\}$ of elements of $P$, we see that the expectation of every finite sum $\sum_{k=1}^nx_k$ is nonnegative and, therefore, the lattice property holds on a set  of positive measure. Thus, we have the induced measure on~$R^P$, which  coincides with $\mu_{\Lambda}$  up to normalization.

2. The restriction of a Bernoulli measure $\mu_{\Lambda}$ to $R^P$ (which is a set of positive measure) is an ergodic central measure. Indeed, as noted above, the restriction of a central measure on $P^{\infty}$, such as the measure $\mu_{\Lambda}$, to $R^P$ is also ergodic and central (in the sense of  $R_P$). Now we check that there are no other such measures. For this, we use the fact that for almost every (in the sense of both measures) point $x$ of $R^P$, the class of~$x$ in the de Finetti partition in the sense of $R^P$  differs from the class of~$x$ in the de Finetti partition in the sense of $P^{\infty}$ by a finite set. One may say that the larger class is fibered over the smaller one with finite fibers.
Further, the conditional measure on the set of additional points is uniform. Therefore, a~central measure on~$R^P$, regarded as a measure on the base, can be uniquely lifted to a~central measure on~$P^{\infty}$, which is ergodic and hence must be a Bernoulli measure. That is, we have established a bijection between the ergodic central measures on the spaces $R^P$ and $P^{\infty}$.
\end{proof}

It follows from the considerations in the previous section that Theorem~1 is a direct corollary of Theorem~2.

As we have indicated, Theorems~1 and~2 remain valid also for the case of infinite posets $P$ and for the case of multiple frequencies. In both cases, the set of lattice sequences does not necessarily have positive Bernoulli measure. Moreover, already in the simplest case of two equal frequencies  $1/2,1/2$, the space $R^P$ is of zero Bernoulli measure.  However, a more careful analysis shows that the bijection from the proof does nevertheless exist.

The trick used in the proof is, of course, quite general. Namely, the centrality property, first, is inherited under restrictions to sets of positive measure and, second, is preserved under extensions of  measures to a covering space provided that the fibers are finite.

The remarkable Thoma theorem describing the characters of the infinite symmetric group can be restated as a theorem describing the central measures on paths in the Young graph. The known proofs of the latter theorem rely on counting tableaux and are related to the theory of symmetric functions. The above proof does not use approximation, but is rather of combinatorial and probabilistic nature. Hence it can be extended to other graphs, for example, continuous graded graphs of Gelfand--Tsetlin type (see \cite{VP,VP1}).

\section{Comments}

\subsection{The problem of describing central measures from a general point of view}

The problem of describing central measures is a special case of the general problem of finding measures with given conditional or pseudo-conditional measures. Stated in another form, it is the problem of describing Markov measures (i.e., their transition probabilities) with given {\it cotransition probabilities}. For more details, see~\cite{V1}. In the general case we have a measured space endowed with a~hyperfinite partition and a cocycle (i.e., a conditional $\sigma$-finite measure),  and the problem is to find all measures with this cocycle (in particular, if the cocycle is identically~$1$, these are the central measures). The specific nature of the  example under consideration is that the hyperfinite partition for one-dimensional central measures is very simple (the de Finetti partition, or the partition into the orbits of the action of the symmetric group), which simplifies things.

\subsection{Relation to Brownian meanders}

The simplest case (where $P$ is a poset consisting of two elements) resembles the  construction of meanders well  known in probability. In this case, $P=\{-1,+1\}$ with $+1>-1$; a sequence $\{x_n\}$ is a lattice sequence if the number of $+1$'s in every initial segment  is greater that the number of $-1$'s, i.e., the corresponding random walk stays in the positive half-plane. If the probability of $+1$ is greater or equal to $\frac12$, then a corresponding central measure is unique and can be interpreted as the measure corresponding to a central random walk in a two-dimensional Weyl chamber.

If $P$ is a chain, then a random numbering corresponds to a random walk in a multidimensional Weyl chamber.

The case $\lambda=1/2$ corresponds to a discrete analog of the notion of {\it Brownian meanders}. However, it is more natural to look for a continuous (diffusion) analog of the notion of central measure; at the moment, such an analog is not known. The case of equal probabilities for a poset with more than two elements is of special interest. The transition probabilities of the corresponding Markov process can be expressed in terms of Schur functions.

\subsection{The Plancherel measure and the GUE (Gaussian unitary ensemble)}

Passing to the limit as the number of elements of the poset $P$ grows in the case of equal probabilities is impossible, it makes sense only if we consider Young diagrams with growing number of rows and columns, and then in the limit we obtain the {\it Plancherel measure}. In fact, such a limiting procedure was performed in our papers with S.~V.~Kerov; the limiting Markov process is described in detail in~\cite{K}. But the intriguing question of how this limiting procedure can be interpreted starting from the finite case and not from the ready Plancherel measure was left undecided. Moreover, it must be of the same nature as passing from unitarily invariant measures on positive definite Hermitian matrices (Wishart measures) to the GUE. The parallelism is that the interplay between  rows and columns
 leads to the Plancherel measure in a similar way  to how the interplay between positive and negative definite Hermitian matrices leads to the GUE. We hope to obtain a new description of the Plancherel and GUE measures, including their limit shape and other characteristics.

\subsection{Multidimensional central measures}

The above comments are directly related to the contents of this paper. However, while in the paper we describe all {\it one-dimensional} central measures, they apply to the description of {\it two-dimensional} central measures. In the cases of the Young graph and the continuous Gelfand--Tsetlin graph (see~\cite{VP}), the two-dimensional invariant measures are, respectively, the Plancherel measure and the GUE. The method described above does not allow one to obtain a direct description of these measures. In the case of~${\Bbb Z}_+^d$ with $d>2$, an even more complicated problem arises of describing central measures of dimensions $3, \ldots, d$. In particular, there is a very long-standing and an extremely important problem concerning the uniqueness of a $d$-dimensional central measure for ${\Bbb Z}_+^d$ (analog of the Plancherel measure for $d=2$) and its analogs for tensors of arbitrary rank.

\subsection{Bernoulli measures and the RSK algorithm}
In \cite {VK}, all central measures for the Young graph were obtained as projections of Bernoulli measures under the Robinson--Schensted--Knuth (RSK) algorithm. In this paper and in~\cite {VP}, is it observed that one-dimensional central measures can be obtained as restrictions of Bernoulli measures to subsets or almost subsets (in the case of equal frequencies). In the theory of invariant measures restrictions are understood as restrictions to invariant subsets, but here this is not the case.
It is possibly because of this distinction that the fact itself has not been observed for a long time. But then a question arises: How the RSK algorithm can be interpreted as a projection to the set of lattice sequences? Is it not an abstract definition of it? Up to now, it was defined only as an algorithm.

\end{document}